\documentclass[10pt,leqno]{amsart}
\usepackage[T1]{fontenc}
\usepackage[utf8]{inputenc}
\usepackage{lmodern}
\usepackage[letterpaper,left=1in,right=1in,top=0.95in,bottom=0.95in]{geometry}
\usepackage{amsmath,amssymb,amsthm,booktabs,array,enumitem,microtype}
\usepackage{hyperref}
\hypersetup{hidelinks,
 pdftitle={A projective Calabi--Yau threefold in characteristic 5 non-liftable to characteristic 0},
 pdfauthor={Dominik Burek, Tymoteusz Chmiel, Slawomir Cynk, Marcin Oczko}}
\setlist{topsep=3pt,itemsep=1pt,leftmargin=18pt}
\newtheorem{theorem}{Theorem}[section]
\newtheorem{proposition}[theorem]{Proposition}
\newtheorem{lemma}[theorem]{Lemma}
\newtheorem{corollary}[theorem]{Corollary}
\theoremstyle{remark}
\newtheorem{remark}[theorem]{Remark}
\title{A projective Calabi--Yau threefold in characteristic 5\\
non-liftable to characteristic 0}
\author{Dominik Burek, Tymoteusz Chmiel, Sławomir Cynk, and Marcin Oczko}
\address{Faculty of Mathematics and Computer Science, Jagiellonian University,
Łojasiewicza 6, 30-348 Kraków, Poland}
\email{dominik.burek@uj.edu.pl}
\email{tymoteusz.chmiel@uj.edu.pl}
\email{slawomir.cynk@uj.edu.pl}
\address{AGH University of Krakow, Faculty of Applied Mathematics
	A. Mickiewicza 30, 30-059 Krakow, Poland}
	\email{marcin.oczko@agh.edu.pl}
\date{}
\begin{document}

\begin{abstract}
We construct a smooth projective Calabi--Yau threefold \(X_5\) over
\(\mathbb F_5\) with \(b_3(X_5)=0\) and \(b_2(X_5)=31\).
Consequently, \(X_5\) has no smooth proper lifting to characteristic zero,
including over ramified bases. To our knowledge, this is the first smooth
projective Calabi--Yau threefold with this property in characteristic
\(p>3\). The threefold is a crepant resolution of the Kummer quotient of
the self-fiber product of a rational elliptic surface. We prove the
vanishing of its third Betti number by semistable reduction and point
counts. A birational double-octic model gives a direct obstruction to
lifting to the ring of Witt vectors of length two, through an explicit
inconsistent system of linear equations for the branch data.
We also determine the Hodge and crystalline cohomology: the threefold is
infinitesimally rigid, its crystalline cohomology is torsion-free, and
its third de Rham cohomology vanishes. Its Hodge--de Rham spectral sequence
has exactly two nonzero differentials, both of rank one on the second page.
\end{abstract}

\maketitle

\section{Introduction}

The existence of smooth projective Calabi--Yau threefolds that do not lift
to characteristic zero is a distinctive feature of positive
characteristic. Hirokado constructed such a threefold in characteristic
three \cite{ref13}, and Schr\"oer gave examples in characteristics
two and three \cite{ref21}. In higher characteristic, small
resolutions produce non-liftable Calabi--Yau algebraic spaces
\cite{ref7}. Obtaining a smooth projective example requires
additional control of the resolution.

In this paper we construct a smooth projective Calabi--Yau threefold
\(X_5\) over \(\mathbb F_5\) with vanishing third Betti number. This rules
out every smooth proper lifting to characteristic zero. To our knowledge,
it is the first smooth projective Calabi--Yau threefold with this
property in characteristic \(p>3\).

It is useful to distinguish this statement from non-liftability to
\(W_2\). The smooth projective characteristic-five example of Cynk and
van Straten has no \(W_2\)-lifting, but it does lift to characteristic zero
over a ramified base \cite[\S5.2]{ref7}.
Our threefold has neither kind of lifting. Since \(5>\dim X_5=3\), it
also lies in the dimension range of the Deligne--Illusie theorem
\cite{ref9}.

Let \(S_+\to\mathbb P^1_{\mathbb F_5}\) be the relatively minimal rational
elliptic surface whose affine Weierstrass equation is
\begin{equation}\label{eq:surface}
 y^2=4x^3+3t(t+2)^3x+3t^4(t+2)^2.
\end{equation}
Its singular fibers have types \(\mathrm{III}\), \(\mathrm{IV}\) and
\(\mathrm{I}_5\), at \(0\), \(3\) and \(\infty\), respectively. Put
\[
 Y_+=S_+\times_{\mathbb P^1}S_+,\qquad
 \iota=(-1,-1),\qquad X_+=Y_+/\langle\iota\rangle.
\]
The existence of a smooth projective crepant resolution of \(X_+\)
follows from the construction of Burek and Cynk
\cite[Theorem~7.4, Remark~7.5 and Example~9.1]{ref4}.
Our main results concern its cohomology and liftability.

\begin{theorem}\label{thm:main}
The quotient \(X_+\) has a smooth projective crepant resolution \(X_5\)
over \(\mathbb F_5\), which is a Calabi--Yau threefold. Its geometric
Betti numbers satisfy
\[
 b_3(X_5)=0,\qquad b_2(X_5)=b_4(X_5)=31.
\]
In particular, \(X_5\) has no smooth proper lifting to characteristic zero.
Put \(k=\overline{\mathbb F}_5\) and
\(X=X_5\times_{\mathbb F_5}k\). Then
\[
 H^1(X,T_{X/k})=0,\qquad H^3_{\mathrm{dR}}(X/k)=0.
\]
The threefold \(X\) has no lifting to \(W_2(k)\), its crystalline
cohomology is torsion-free, and Frobenius acts as zero on
\(H^3(X,\mathcal O_X)\).
\end{theorem}

Thus \(X_5\) provides a characteristic-five example of the type sought
in Schoen's Question~10.2(i) \cite{ref20}.

The proof uses two complementary descriptions of the threefold.
Sections~3--5 study the self-fiber product and its Kummer quotient.
A semistable model and two finite-field point counts give the vanishing
of the third Betti number of a resolved self-product. The fixed curves
of diagonal negation are rational, so resolving the quotient introduces
no third cohomology.

Section~6 presents a birational double-octic model. We first exclude a
\(W_2\)-lifting of a projective small resolution by an explicit rank
calculation for the lifted branch data. A second calculation proves its
infinitesimal rigidity. The obstruction class and the conjugate spectral
sequence then give the vanishing of third de Rham cohomology, as in
Lam's criterion \cite[Corollary~1.6 and Lemmas~5.1--5.2]{ref16}, using the
description of the differential due to Achinger and Suh
\cite[Theorem~1.3 and Corollary~4.4]{ref1}.
Integral crystalline correspondences transfer the middle cohomology
calculation to \(X_5\), and the Deligne--Illusie theorem excludes
\(W_2\)-liftings of this model as well.
Section~7 derives the remaining Hodge, crystalline and deformation
consequences. In particular, we determine the two nonzero differentials
in the Hodge--de Rham spectral sequence.

\par\medskip\noindent\textbf{Notation and conventions.}
A Calabi--Yau threefold over a field is a smooth proper geometrically
connected threefold with trivial canonical bundle and
\(H^1(\mathcal O)=H^2(\mathcal O)=0\).
Projectivity is stated separately when algebraic spaces occur.
Betti numbers are geometric \(\ell\)-adic Betti numbers, with \(\ell\)
different from the characteristic.
The ground field is specified locally: the computations in Section~6
are over \(\mathbb F_5\), whereas Section~7 is over
\(\overline{\mathbb F}_5\), with the indicated base changes understood.
We write \(W(k)\) for Witt vectors and
\(W_2(k)=W(k)/(25)\) when \(\operatorname{char}(k)=5\).
In the mixed-characteristic model, a subscript \(\infty\) denotes the
generic fiber and a subscript \(+\) the special fiber; the base point
\(t=\infty\) is distinguished by the variable \(t\).
A superscript \(\mathrm W\) denotes a Weierstrass model.

\subsection*{Acknowledgements}
Tymoteusz Chmiel was supported by the National Science Center, Poland,  grant no. 2023/49/N/ST1/04089.
S\l awomir Cynk was supported by National Science Centre, Poland, grant No. 2020/39/B/ST1/03358.
Marcin Oczko was supported by the National Science Centre, Poland,
grant No.~2025/57/N/ST1/05021.

The authors used OpenAI’s ChatGPT and Codex to help check calculations and polish the exposition. The authors take responsibility for the mathematical claims and the final text.

\section{Preliminaries}

\subsection{Birational Calabi--Yau threefolds over finite fields}

\begin{proposition}\label{prop:codim-one}
Let $k$ be a perfect field, and let $X$ and $X'$ be smooth, proper, integral algebraic spaces of finite type
  over $k$ with trivial canonical bundles. Every birational map $\varphi\colon X\dashrightarrow X'$ is an
  isomorphism in codimension one: there are closed subspaces $Z\subset X$ and $Z'\subset X'$, both of
  codimension at least two, such that $\varphi$ restricts to an isomorphism $X\setminus Z\simeq X'\setminus
  Z'$.
\end{proposition}

\begin{proof}
Identify the function fields using $\varphi$, choose nowhere-vanishing generators $\eta$ and $\eta'$ of the
  canonical bundles, and write $\eta'=f\eta$ as rational top forms. Properness of $X'$ implies that $\varphi$
  is defined at every codimension-one point of $X$, so $f$ has no divisorial pole. Normality makes $f$
  regular, and properness and integrality imply that this nonzero function is a unit.

If a prime divisor $D\subset X$ were contracted, its center on the smooth space $X'$ would have codimension
  at least two. Its discrepancy over $X'$ would then be strictly positive, giving
  $\operatorname{ord}_D(\eta')>0$. This contradicts $\eta'=f\eta$, since $f$ is a unit and $\eta$ is nowhere
  vanishing. Applying the same argument to $\varphi^{-1}$ shows that neither map contracts a prime divisor.
  They are therefore inverse isomorphisms on open subspaces whose complements have codimension at least two.
\end{proof}

\begin{corollary}\label{cor:betti-birational}
Let $k$ be a finite field, and let $X$ and $X'$ be smooth proper integral three-dimensional algebraic spaces
  over $k$ with trivial canonical bundles. If $X$ and $X'$ are birational, they have equal Betti numbers.
\end{corollary}

\begin{proof}
By Proposition~\ref{prop:codim-one}, there are isomorphic open subspaces whose complements $Z$ and $Z'$ have
  dimension at most one. The localization sequences give the equality of virtual Galois representations
\[
[H^*(X)]-[H^*(X')]=[H^*(Z)]-[H^*(Z')],
\]
where brackets denote alternating sums. All weights on the right are at most two. The cohomology of a smooth
  proper algebraic space is pure, so the multisets of Frobenius eigenvalues on $H^i(X)$ and $H^i(X')$ agree
  for $i\ge3$. Purity for algebraic spaces follows by proper cohomological descent and Poincar\'e duality
  from the scheme case of \cite{ref8}. Poincar\'e duality then gives equality of the Betti numbers in degrees
  zero, one and two as well.
\end{proof}

\subsection{Cohomology of the resolved involution quotient}

\begin{lemma}\label{lem:quotient-cohomology}
Let $X$ be a smooth proper three-dimensional algebraic space over a finite field of odd characteristic, and
  let $\iota$ be an involution with fixed locus $C$ of pure dimension one. Let $Y$ be the blow-up of the
  image of $C$ in $X/\langle\iota\rangle$. For $\ell$ different from the characteristic, there is a
  Galois-equivariant isomorphism
\[
H^m(Y_{\bar k},\mathbb Q_\ell)
\simeq H^m(X_{\bar k},\mathbb Q_\ell)^\iota
\oplus H^{m-2}(C_{\bar k},\mathbb Q_\ell)(-1).
\]
\end{lemma}

\begin{proof}
The fixed curve is smooth, with normal action $(-1,-1)$. In completed local coordinates $x,y$ normal to $C$,
  the quotient has invariants $a=x^2$, $b=xy$, $c=y^2$ with $ac=b^2$. On the chart $y=xv$ of
  $\operatorname{Bl}_C X$ the invariants are $a=x^2,v$, with $b=av$, $c=av^2$. This is the $a$-chart of the
  blow-up of $(a,b,c)$; the chart with $c$ is symmetric. Thus $Y=(\operatorname{Bl}_C
  X)/\langle\iota\rangle$. For this finite quotient the invariant direct image of the constant $\mathbb
  Q_\ell$-sheaf is the constant sheaf, and its higher direct images vanish. Taking invariants is exact. The
  blow-up formula then gives the displayed decomposition, since the involution fixes $C$ and acts trivially
  on the projectivization of its normal bundle.
\end{proof}

\section{Fiber product}

Consider the elliptic surface $\mathcal S^W$ over $\operatorname{Spec}\mathbb Z_5$ defined by the Weierstrass
  equation
\[
 y^2=4x^3+A(t)x+B(t)
\]
where
\[
\begin{aligned}
 A(t)&:=-\frac{3}{64}(t^2-28t+100)(t+2)^2,\\
 B(t)&:=\frac{1}{512}(t^4-40t^3+360t^2+160t+2000)(t+2)^2.
\end{aligned}
\]
The generic Weierstrass surface has rational double points of types $A_1$, $A_2$ and $A_4$ above $0$, $-2$
  and $\infty$, respectively. The minimal regular elliptic surfaces of the generic and special Weierstrass
  models are denoted by $S_\infty$ and $S_+$.

Lemma~\ref{lem:partial-resolution} constructs a projective simultaneous minimal resolution $\mathcal S$ over
  $\mathbb Z_5$, with fibers $S_\infty$ and $S_+$.

The discriminant in this Weierstrass normalization is
\[
 \Delta_{\mathcal S^W}=-A(t)^3-27B(t)^2
 =\frac{3^6}{2^6}\,t^2(t-25)(t+2)^4.
\]
By the Tate algorithm, the surface $S_\infty$ has four singular fibers (see \cite{ref12}):
\[
\begin{array}{c|cccc}
 t&-2&0&25&\infty\\\hline
 S_{\infty,t}&\mathrm{IV}&\mathrm{I}_2&\mathrm{I}_1&\mathrm{I}_5.
\end{array}
\]
Modulo $5$ the singular fibers of types $\mathrm{I}_1$ and $\mathrm{I}_2$ collide and the surface $S_+$ has
  three singular fibers
\[
\begin{array}{c|ccc}
 t&3&0&\infty\\\hline
 S_{+,t}&\mathrm{IV}&\mathrm{III}&\mathrm{I}_5.
\end{array}
\]
We consider the two self-fiber products
\[
 Y_\infty=S_\infty\times_{\mathbb P^1}S_\infty,
 \qquad Y_+=S_+\times_{\mathbb P^1}S_+.
\]
The characteristic-five surface is the one in \eqref{eq:surface}; compare \cite[Example~6.1]{ref4}.

\begin{proposition}\label{prop:generic}
After a finite extension of $\mathbb Q_5$, the generic fiber $Y_\infty$ admits a projective crepant
  resolution $\widetilde Y_\infty$ with geometric Hodge numbers
\[
 h^{1,1}(\widetilde Y_\infty)=42,
 \qquad h^{1,2}(\widetilde Y_\infty)=0
\]
and the third Betti number
\[
 b_3(\widetilde Y_\infty)=2.
\]
\end{proposition}
\begin{proof}
The four common singular fibers give the local products
\[
 \mathrm{IV}\times\mathrm{IV},\qquad
 \mathrm{I}_2\times\mathrm{I}_2,\qquad
 \mathrm{I}_1\times\mathrm{I}_1,\qquad
 \mathrm{I}_5\times\mathrm{I}_5.
\]
The product $\mathrm{IV}\times\mathrm{IV}$ contains an isolated ordinary triple point and admits a projective
  crepant resolution by the blow-up of the singular point. For $\mathrm{I}_n\times\mathrm{I}_n$, with $n>1$,
  projective small crepant resolutions can be realized by consecutive blow-ups of components of the singular
  fiber. In the self-product case the node arising from $\mathrm{I}_1\times\mathrm{I}_1$ also admits a
  projective small resolution, obtained globally by blowing up the diagonal. Consequently, $Y_\infty$, after
  a finite extension of its ground field, admits a smooth projective crepant resolution
\[
 \widetilde Y_\infty\longrightarrow Y_\infty.
\]
For a resolved common fiber, denote by $N_c$ its number of irreducible components and by $e$ the Euler
  characteristic. The local contributions are (\cite[Theorem~4.1]{ref4})
\[
\begin{array}{c|rr}
 \text{fiber product}&N_c&e\\\hline
 \mathrm{IV}\times\mathrm{IV}&10&24\\
 \mathrm{I}_2\times\mathrm{I}_2&4&8\\
 \mathrm{I}_1\times\mathrm{I}_1&1&2\\
 \mathrm{I}_5\times\mathrm{I}_5&25&50.
\end{array}
\]
The table gives $40$ components in the resolved common fibers and Euler characteristic $84$. In the
  Hodge-number formula \cite[Proposition~5.2]{ref4}, the union of the singular-fiber sets has cardinality
  $4$, each surface has one additive fiber, and there are no fibers singular on just one factor. The
  endomorphism rank of the generic elliptic curve is $1$ by non-isotriviality in characteristic zero. Hence
\[
 h^{1,1}=1+4-5+1+1+40=42,
 \qquad h^{1,2}=42-\frac{84}{2}=0.
\]
The Calabi--Yau Hodge numbers give $b_3=2$.
\end{proof}

\section{The fiber product in characteristic 5}

Put $Y_+=S_+\times_{\mathbb P^1_{\mathbb F_5}}S_+$ and let $Y_+^W$ be the corresponding product of
  Weierstrass models. The common singular fibers have types $\mathrm{III}\times\mathrm{III}$,
  $\mathrm{IV}\times\mathrm{IV}$ and $\mathrm{I}_5\times\mathrm{I}_5$. We construct a projective partial
  resolution in mixed characteristic and its semistable reduction. On the special fiber we use both a crepant
  small resolution in algebraic spaces and a projective big resolution with one exceptional quadric. The big
  resolution is not crepant.

\begin{lemma}\label{lem:partial-resolution}
Let $R=W(\mathbb F_5)=\mathbb Z_5$, with fraction field $K=\mathbb Q_5$. The Weierstrass surface $\mathcal
  S^W$ admits a projective simultaneous minimal resolution $\mathcal S\to\mathcal S^W$, smooth over $R$.
  There is a projective modification
\[
 \mathcal P\longrightarrow\mathcal Y:=\mathcal S\times_{\mathbb P^1_R}\mathcal S
\]
whose generic fiber is a smooth projective crepant resolution of $Y_\infty\otimes K$. Its special fiber $P$
  is a crepant partial resolution of $Y_+$ with precisely one ordinary double point $N$. The scheme $\mathcal
  P$ is regular, is smooth over $R$ away from $N$, and at $N$ its completed local equation is
\begin{equation}\label{eq:node}
 cd=u^2+v^2-2a,\qquad a\in R,\qquad v_5(a)=1.
\end{equation}
In particular $B=\operatorname{Bl}_N P$ is a smooth projective threefold over $\mathbb F_5$.
\end{lemma}

\begin{proof}
Put $R_{\mathrm{ur}}=W(\mathbb F_{25})$; we use this extension only to check split charts. The singular
  sections of $\mathcal S^W$ over $t=0,-2,\infty$ are disjoint. At $t=0$, write $x=X+5/4$ and blow up
  $(X,y,t)=0$. The tangent conic
\[
 y^2=15X^2-\frac{27}{2}Xt+\frac{135}{16}t^2
\]
is nondegenerate over $R$, so this resolves the relative $A_1$ singularity. At $t=-2$, with $s=t+2$, the
  equation is $y^2=4x^3+s^2(\alpha(s)x+\beta(s))$, where $\beta(0)=27/4$. Over $R_{\mathrm{ur}}$, the
  substitutions $U=y-s\sqrt{\alpha x+\beta}$ and $V=y+s\sqrt{\alpha x+\beta}$ give $UV=4x^3$. Blowing up
  $(x,y,s)=0$ resolves this relative $A_2$ singularity; its exceptional components meet in a unique section
  defined over $R$.

At infinity, put $s=1/t$. The singular section is $(x,y,s)=(1/16,0,0)$, the quadratic coefficient at the
  double root is $3/4$, and the discriminant has order five. Completing squares over $R_{\mathrm{ur}}$ gives
  $UV=s^5$. Blow up this section and then the remaining singular section: the successive $s$-charts are
  $U_1V_1=s^3$ and $U_2V_2=s$, and the other charts are smooth. The second center is the intersection section
  of the first exceptional components and descends to $R$. These global blow-ups give the simultaneous
  minimal resolution $\mathcal S$. All surface centers are invariant under negation, which therefore lifts to
  $\mathcal S$.

Over $t=-2$, blow up the product of the exceptional intersection section with itself. On the surface chart
  $U=xr$, $V=xw$, one has $x=rw/4$ and $s=rw(w-r)\cdot\text{unit}$. Thus the product has an ordinary triple
  point along the section, with tangent cubic $rw(w-r)-r'w'(w'-r')$. Its projectivization is smooth, since
  $2,3$ are units and the two derivatives of each binary cubic vanish simultaneously only at the origin. The
  blow-up resolves the product there.

At infinity, label the split five-cycle $0,1,2,2',1'$. Galois and negation act by the same reflection. Omit
  $(0,2),(0,2'),(2,0),(2',0)$ and blow up the strict transforms of the unions in each remaining reflection
  orbit, as in \cite[Theorem~7.4, Remark~7.5]{ref4}. These invariant unions descend to $R$. Every node lies
  on a retained component. Distinct components in one orbit are either disjoint or opposite at a node: the
  only adjacent reflected vertices are $2,2'$, and the omissions exclude a fixed other coordinate. At
  $xy=zw$, the opposite ideals satisfy
\[
 I=(x,z),\quad J=(y,w),\qquad I\cap J=IJ,\qquad J=(w/x)I.
\]
The union ideal is a fractional principal multiple of $I^2$, so its blow-up is the usual small resolution.
  Subsequent strict transforms are Cartier near that node. This resolves all twenty-five node sections
  equivariantly over $R$.

Finally blow up the strict diagonal and the strict graph of negation. These are global Weil divisors, Cartier
  on the locus already smooth over $R$. For the collision over $t=0$, use the surface chart $y=Xu$, $t=Xv$.
  Its equation is $u^2=G(X,v)$, where
\[
 G(X,0)=4X+15,\qquad G_v(0,0)=-27/2,\qquad
 G\equiv4X+4v\pmod{(5,(X,v)^2)}.
\]
Thus $v=F(X,s)$, with $s=u^2$ and $F\equiv4s-X\pmod{(5,(X,s)^2)}$. The parameter Morse lemma gives
  $t=z^2+h(s)$: the quadratic coefficient has square residue $4$. This coordinate change fixes $s$ and
  commutes with negation. Since $F(0,15)=0$ and $F_X(0,15),F_s(0,15)$ are units,
\[
 h(15)=h'(15)=0,\qquad
 h''(15)=-\frac{F_s(0,15)^2}{2F_X(0,15)}.
\]
Evaluation and translation by $15$ are legitimate in the complete ring $R[[s]]$. Hence $h(s)=(s-15)^2H(s)$,
  with $H(0)\equiv4\pmod5$. Choose $c(s)^2=-H(s)$ with $c(0)\equiv1\pmod5$, and put $g(s)=(s-15)c(s)$ and
  $a=-g(0)$. Then $v_5(a)=1$ and $g'(0)\equiv1\pmod5$. Replacing $u$ by $u\sqrt{(g(u^2)-g(0))/u^2}$ gives,
  compatibly with negation,
\begin{equation}
 t=z^2-(u^2-a)^2.\label{eq:collision}
\end{equation}
For the self-product, put $A=z_1-z_2$ and $d=z_1+z_2$. Its equation is
\[
 Ad=(u-v)(u+v)(u^2+v^2-2a).
\]
Blowing up the diagonal $(A,u-v)$ gives a smooth chart and the residual chart
\[
 A=(u-v)r,\qquad rd=(u+v)(u^2+v^2-2a).
\]
The strict antidiagonal has ideal $(r,u+v)$. Its blow-up again has a smooth chart and, on $r=(u+v)c$, the
  chart $cd=u^2+v^2-2a$ of (2). The chart rings are $R$-flat and reduce to the asserted partial-resolution
  charts, identifying the special fiber directly.

The special fiber of \eqref{eq:node} has one ordinary double point, while the total space is regular because
  $v_5(a)=1$. Its generic fiber is smooth. The critical points $z=0$, $u=0,\pm\sqrt a$ in
  \eqref{eq:collision} account for the four generic nodes over $0$ and the node over $25$; the other surface
  charts are smooth for the elliptic map away from its critical locus. The centers over $-2$ and infinity are
  disjoint from this neighborhood. All centers are closed algebraic subschemes, so the construction is
  projective. The triple-point blow-up and the small resolutions are crepant on both fibers. Blowing up the
  last node of $P$ replaces it by its smooth tangent quadric and gives $B$.
\end{proof}

\begin{lemma}\label{lem:semistable}
Let $R'=\mathbb Z_5[\pi]$, where $\pi^2=5$. After base change to $R'$, the generic fiber of $\mathcal P$ has
  a projective strictly semistable model $\mathcal Z/R'$ with special fiber
\[
 \mathcal Z_+=B\cup Q,\qquad B\cap Q=V.
\]
Here $Q$ is a smooth quadric threefold and $V$ is a smooth quadric surface. The component $B$ is the
  resolution in Lemma~\ref{lem:partial-resolution}.
\end{lemma}
\begin{proof}
Blow up the closed point $N$ of the regular four-dimensional scheme $\mathcal P$. The exceptional divisor is
  $\mathbb P^3$, and the special fiber is
\[
 B+2\mathbb P^3.
\]
Indeed, by (2) the uniformizer has order two at $N$, with nondegenerate quadratic initial form. Its zero
  locus on the exceptional $\mathbb P^3$ is $V$, and the strict transform $B$ meets this divisor transversely
  along $V$.

Define $\mathcal Z$ as the normalization of this blowup after the base change to $R'$. Normalization is
  finite here, since the rings are excellent, so the resulting model is projective. Along $V$ the equation
  before base change is $xy^2=5$, after absorbing a unit into the local equation $x$ of $B$. The
  normalization is given by
\[
 x=z^2,\qquad yz=\pi.
\]
This is a strict semistable crossing. Away from $V$ the same calculation gives smooth charts. The component
  lying over $B$ is isomorphic to $B$. The other component is the double cover of $\mathbb P^3$ branched
  along $V$, hence a smooth quadric threefold $Q\subset\mathbb P^4$. This cover is supplied globally by the
  normalized base change, so no separate square root of a divisor line bundle needs to be chosen. Its
  intersection with $B$ is $V$.
\end{proof}

\begin{theorem}
Let $S\to\mathbb P^1$ be the rational elliptic surface over $k=\overline{\mathbb F}_5$ with singular fibers
  $\mathrm{III}$, $\mathrm{IV}$ and $\mathrm{I}_5$ at $0,3,\infty$, respectively, and put
  $Y_+=S\times_{\mathbb P^1}S$. There is a smooth proper crepant resolution $\widehat Y_+\to Y_+$ in
  algebraic spaces for which
\[
 b_1(\widehat Y_+)=b_3(\widehat Y_+)=b_5(\widehat Y_+)=0,
 \qquad b_2(\widehat Y_+)=b_4(\widehat Y_+)=42.
\]
The projective big resolution $B$ has Betti numbers
\[
 (b_0(B),\ldots,b_6(B))=(1,0,43,0,43,0,1).
\]
\end{theorem}
\begin{proof}
Apply the Rapoport--Zink weight spectral sequence to $\mathcal Z/R'$ from Lemma~\ref{lem:semistable}; see
  \cite[Satz~2.10]{ref17} and \cite[Corollary~2.2.4, Proposition~2.2.6]{ref18}. All cohomology groups in this
  proof are geometric \'etale cohomology with coefficients in $\mathbb Q_\ell$, $\ell\ne5$. Put $\mathcal
  Z^{(0)}=B\sqcup Q$, $\mathcal Z^{(1)}=V$, and $\mathcal Z^{(j)}=\varnothing$ for $j\ge2$. The sequence is
\begin{equation*}
 E_1^{p,q}=\bigoplus_{i\ge\max(0,-p)}
 H^{q-2i}\bigl(\mathcal Z^{(p+2i)}\bigr)(-i)
 \quad\Longrightarrow\quad H^{p+q}(\mathcal Z_{\overline K}).
\end{equation*}
Its only nonzero columns are
\[
\begin{array}{c|ccc}
 p&-1&0&1\\\hline
 E_1^{p,q}&H^{q-2}(V)(-1)&H^q(B)\oplus H^q(Q)&H^q(V).
\end{array}
\]
The first differential consists of the signed Gysin and restriction maps. Each term in row $q$ is pure of
  weight $q$, so Frobenius equivariance and purity imply degeneration at $E_2$. By
  Proposition~\ref{prop:generic}, the geometric Betti numbers of the generic fiber are $(1,0,42,2,42,0,1)$.

The map $H^0(B)\oplus H^0(Q)\to H^0(V)$ is surjective, and $H^0(V)(-1)\to H^2(B)\oplus H^2(Q)$ is injective:
  its component in $H^2(Q)$ is the nonzero divisor class of $V$. Since $Q$ and $V$ have no odd cohomology,
  degree one of (4) gives $H^1(B)=0$. Poincar\'e duality then gives $H^5(B)=0$.

Let
\[
 r=\operatorname{rank}\bigl(H^2(B)\oplus H^2(Q)\longrightarrow H^2(V)\bigr),
 \qquad \alpha=2-r.
\]
The hyperplane class of $Q$ restricts nontrivially to $V$, so $r\in\{1,2\}$. Duality identifies $r$ with the
  rank of the Gysin map $H^2(V)(-1)\to H^4(B)\oplus H^4(Q)$. In total degrees two and three the $E_2$ terms
  therefore give
\begin{equation}\label{eq:betti-rel}
 b_2(B)+\alpha-2=42,\qquad b_3(B)+2\alpha=2.
\end{equation}
Thus $(b_2(B),b_3(B))$ is either $(43,0)$ or $(44,2)$.

\smallskip\noindent\emph{Point counts.}
Write $q=5^r$, with $r\ge1$. Only the counts for $q=25,625$ enter the argument. Let $\chi_q$ be the quadratic
  character of $\mathbb F_q$, with $\chi_q(0)=0$. For $t\in\mathbb F_q$, set
\[
 n_t=q+1+\sum_{x\in\mathbb F_q}
 \chi_q\bigl(4x^3+3t(t+2)^3x+3t^4(t+2)^2\bigr).
\]
At infinity the equation is $y^2=4x^3+3x+3$, giving $n_\infty$ by the same rule. Direct evaluation of
  $\sum_{t\in\mathbb P^1(\mathbb F_q)}n_t^2$ in Magma gives $19401$ and $246485001$ for $q=25,625$,
  respectively.

Write $\epsilon_r=(-1)^r$. Resolving the Weierstrass surfaces changes only the three singular fibers. Their
  point counts are
\[
\begin{array}{c|cc}
 t&\#S_{+,t}^W(\mathbb F_q)&\#S_{+,t}(\mathbb F_q)\\\hline
 0&q+1&2q+1\\
 3&q+1&(2+\epsilon_r)q+1\\
 \infty&q+1-\epsilon_r&(3+2\epsilon_r)q+1-\epsilon_r.
\end{array}
\]
At infinity, $4x^3+3x+3=4(x-1)^3+2(x-1)^2$; the two tangents split over $\mathbb F_{25}$. The two nonidentity
  components of the type $\mathrm{IV}$ fiber also split over $\mathbb F_{25}$. Taking the differences of the
  squares of the entries gives
\[
 \#Y_+(\mathbb F_q)-\#Y_+^W(\mathbb F_q)
 =(19+16\epsilon_r)q^2+(4+2\epsilon_r)q,
\]
which gives $\#Y_+(\mathbb F_{25})=41426$ and $\#Y_+(\mathbb F_{625})=260160626$.

Over $\mathbb F_{25}$ the ordinary triple point is replaced by the smooth cubic surface
\[
 E:\quad x(x^2-2y^2)-z(z^2-2w^2)=0\subset\mathbb P^3.
\]
To count it directly, set $m_b=\#\{(x,y)\in\mathbb F_q^2:x^3-2xy^2=b\}$. Then
\[
 \#E(\mathbb F_q)=\frac{\sum_{b\in\mathbb F_q}m_b^2-1}{q-1}.
\]
Evaluation for $q=25,625$ gives $801$ and $395001$. All twenty-five nodes over infinity are split and are
  replaced by projective lines. The diagonal and antidiagonal each add one projective line over the
  collision. The final node is replaced by a split quadric surface: its equation is $cd=u^2+v^2$, and $-1$ is
  a square in $\mathbb F_{25}$. The point-count increments are therefore
\[
\begin{array}{r|rrrrr}
 q&\#Y_+(\mathbb F_q)&\#E(\mathbb F_q)-1&25q&2q&q^2+2q\\\hline
 25&41426&800&625&50&675\\
 625&260160626&395000&15625&1250&391875.
\end{array}
\]
Adding the increments gives $\#B(\mathbb F_{25})=43576$ and $\#B(\mathbb F_{625})=260964376$. For both fields
  this is $q^3+43(q^2+q)+1$.

\smallskip\noindent\emph{Excluding $b_3(B)=2$.}
For $q=25,625$, let $t_i(q)$ be the trace of geometric $q$-Frobenius on $H^i(B)$. Hard Lefschetz for an ample
  class defined over $\mathbb F_{25}$ gives $t_4(q)=qt_2(q)$, and $t_2(q)$ is an integer by the Weil
  conjectures \cite{ref8}. The trace formula and the two point counts yield
\[
 t_3(q)=(q+1)\bigl(t_2(q)-43q\bigr).
\]
If $b_3(B)=2$, write $t_3(25)=26m$ with $m\in\mathbb Z$. The Weil bound gives $|26m|\le2\cdot25^{3/2}=250$,
  hence $|m|\le9$. The alternating Poincar\'e pairing on this two-dimensional space has Frobenius multiplier
  $25^3$. Therefore the determinant of Frobenius is exactly $25^3$, and
\[
 t_3(625)=t_3(25)^2-2\cdot25^3=676m^2-31250.
\]
The point count at $625$ makes this integer divisible by $626$. Reduction modulo $313$ gives $313\mid m^2+1$,
  which is impossible since $1\le m^2+1\le82$. Consequently $b_3(B)=0$, and \eqref{eq:betti-rel} gives
  $b_2(B)=43$.

Either small resolution $\widehat Y_+\to P$ of the last node exists as a proper algebraic space; equivalently
  one applies the usual algebraic-space modification of an ordinary double point \cite{ref2}. It is crepant
  and
\[
 B\simeq\operatorname{Bl}_C\widehat Y_+,
 \qquad C\simeq\mathbb P^1,
 \qquad N_{C/\widehat Y_+}\simeq\mathcal O(-1)^{\oplus2}.
\]
The blowup formula subtracts one from $b_2$ and $b_4$ and leaves the odd Betti numbers unchanged. This proves
  the theorem.
\end{proof}

\begin{corollary}\label{cor:betti-small}
Every smooth proper crepant resolution of $Y_+$ has third Betti number zero.
\end{corollary}
\begin{proof}
The dualizing bundle of $Y_+$ is trivial by the canonical bundle formula for a fiber product of rational
  elliptic surfaces. Thus every smooth crepant resolution has trivial canonical bundle. Apply
  Corollary~\ref{cor:betti-birational}, after a finite extension of the field of definition if necessary, to
  the resolution constructed above.
\end{proof}

\section{Kummer fibration}

Let $\iota$ be diagonal fiberwise negation on $Y_+$, and put $X_+=Y_+/\langle\iota\rangle$. We compare a
  projective crepant model of $X_+$ with the quotient of an equivariant crepant algebraic-space resolution of
  $Y_+$.

\begin{lemma}\label{lem:equivariant}
The self-product $Y_+=S_+\times_{\mathbb P^1}S_+$ admits a smooth proper crepant resolution $\widehat Y\to
  Y_+$ in algebraic spaces to which diagonal negation lifts. It can be chosen over $\mathbb F_5$.
\end{lemma}
\begin{proof}
At $\mathrm{I}_5\times\mathrm{I}_5$, use the special fibers of the paired blowups in
  Lemma~\ref{lem:partial-resolution}. Their split local charts reduce to the same small-resolution charts
  over $\mathbb F_{25}$, and the descended centers are invariant under diagonal negation. The blowup of the
  ordinary triple point at $\mathrm{IV}\times\mathrm{IV}$ is canonical, smooth and crepant, and the
  involution lifts.

At $\mathrm{III}\times\mathrm{III}$, start with
\[
 x(x-u^2)-y(y-v^2)=0,\qquad
 \iota(x,y,u,v)=(x,y,-u,-v).
\]
Put $a=2x-u^2$, $b=2y-v^2$, $X=a-b$ and $Y=a+b$. The equation becomes
\[
 XY=u^4-v^4=\prod_{j=1}^4l_j(u,v),
\]
where the four linear factors are defined over $\mathbb F_5$ and $\iota(X,Y,u,v)=(X,Y,-u,-v)$. Blow up the
  ideal $(X,l_4)$. Its complementary chart is smooth, while the chart $X=p_1l_4$ has equation
\[
 p_1Y=l_1l_2l_3,\qquad
 \iota(p_1,Y,u,v)=(-p_1,Y,-u,-v).
\]
Continue with $(p_1,l_3)$ and then $(p_2,l_2)$, where $p_2=p_1/l_3$ on the residual chart. The sign of the
  first coordinate alternates at each step, so every centre is invariant. Each complementary chart is smooth,
  and the final residual equation is linear in $u,v$. This gives a small crepant resolution with exceptional
  fiber a chain of three projective lines. The construction is defined over $\mathbb F_5$. Artin's
  algebraization of formal modifications \cite[Theorem~3.2]{ref2} gives the required algebraic-space
  modification, and the equivariant local modifications glue with the unchanged smooth locus.
\end{proof}

\begin{lemma}\label{lem:fixed-curves}
Every geometric connected component of $C=\operatorname{Fix}(\iota|_{\widehat Y})$ is a smooth rational
  curve, and there are exactly seven such components.
\end{lemma}
\begin{proof}
The involution preserves the volume form: it changes the signs of both elliptic relative one-forms and fixes
  the base differential. Since $2\ne0$, its fixed locus is smooth. At a fixed point its determinant is one;
  nontriviality excludes a three-dimensional fixed tangent space. The tangent signs are therefore
  $(1,-1,-1)$, so $C$ is purely one-dimensional.

The nonzero two-torsion multisection has equation
\[
 D:\qquad4x^3+3t(t+2)^3x+3t^4(t+2)^2=0.
\]
The off-diagonal part of $D\times_{\mathbb P^1}D$ is given by this equation and
\[
 4(x^2+xy+y^2)+3t(t+2)^3=0.
\]
It has the parametrization
\[
 t=\frac{4(r^2-1)^2}{r^2(r^2+1)^2},\qquad
 x=\frac{3(r^2-1)(r^2-2)^2}{r^4(r^2+1)^4},\qquad
 y=2(r^5+1)x.
\]
Here is a direct verification. Set
\[
\begin{gathered}
 T=4(r^2-1)^2,\qquad d=r^2(r^2+1)^2,\qquad
 X=3(r^2-1)(r^2-2)^2,\\
 Y=2(r^5+1)X,\qquad U=T+2d.
\end{gathered}
\]
In $\mathbb F_5[r]$ one has
\[
 4X^3+3TU^3X+3T^4U^2=0,\qquad
 4(X^2+XY+Y^2)+3TU^3=0.
\]
These are exactly the two equations after clearing denominators. Moreover, $t(r)$ has degree six and
  $r^5=y/(2x)-1$ belongs to $L=k(t,x,y)\subset k(r)$. Consequently $[k(r):L]$ divides both five and six and
  is one. The generic algebra of the two equations has rank six over $k(t)$, so this is the entire
  off-diagonal component, and its normalization is rational. Since $t,x$ depend on $r^2$ and $y$ satisfies a
  quadratic over $k(t,x)$, one also has $k(t,x)=k(r^2)$; hence $D$ is rational. These arguments hold after
  algebraic closure as well.

The horizontal fixed locus therefore has five rational components: the zero-section product, the two products
  of the zero section with $D$, the diagonal $D$, and the off-diagonal component. We must also check the
  vertical components. On each of the fibers $\mathrm{I}_5$, $\mathrm{III}^{b}$ and $\mathrm{IV}$, negation
  is nontrivial on every preserved component and permutes the other components; see \cite[Lemma~7.1]{ref4}.
  Thus every vertical fixed curve on $\widehat Y$ is exceptional. At $\mathrm{IV}\times\mathrm{IV}$, the
  exceptional cubic has equation $x(x^2-au^2)-z(z^2-av^2)=0$, with $a\ne0$. The involution changes the signs
  of $u,v$. The fixed linear subspace $u=v=0$ meets the cubic in three points, while $x=z=0$ gives its fixed
  line. Thus this fiber contributes one rational curve. At $\mathrm{III}\times\mathrm{III}$, the exceptional
  projective lines of the three factor blowups have coordinates $[X:l_4]$, $[p_1:l_3]$ and $[p_2:l_2]$, on
  which the signs are respectively $(+,-)$, $(-,-)$ and $(+,-)$. Exactly the middle component is fixed
  pointwise. At $\mathrm{I}_5\times\mathrm{I}_5$ only one of the twenty-five nodes is fixed. At that node the
  involution exchanges both pairs of branches. Write its equation as $\det(M)=0$, where
\[
 M=\begin{pmatrix}x&u\\v&y\end{pmatrix},\qquad
 \iota(M)=JMJ,\qquad
 J=\begin{pmatrix}0&1\\1&0\end{pmatrix}.
\]
On either invariant small resolution, the exceptional projective line has the nontrivial action induced by
  $J$. The other exceptional curves are permuted, so this fibre contributes no vertical fixed curve. The
  preceding fixed-locus analysis excludes vertical curves away from the exceptional locus. Hence there are
  seven components and $e(C)=14$.
\end{proof}

\begin{proposition}\label{prop:betti-kummer}
The quotient $X_+=Y_+/\langle\iota\rangle$ has a smooth projective Calabi--Yau crepant resolution $X_5$ over
  $\mathbb F_5$. It satisfies $b_3(X_5)=0$ and $b_2(X_5)=b_4(X_5)=31$, and has no smooth proper lifting to
  characteristic zero.
\end{proposition}
\begin{proof}
For projectivity use \cite[Theorem~7.4, Remark~7.5]{ref4}: negation is of type $\mathrm{III}^{b}$, since it
  fixes the zero section and preserves its fiber component. The relevant products are
  $\mathrm{I}_5\times\mathrm{I}_5$, $\mathrm{III}^{b}\times\mathrm{III}^{b}$, and
  $\mathrm{IV}\times\mathrm{IV}$, all in the projective list. One can also identify the precise example over
  the ground field. Over $\mathbb F_{25}$ the substitution
\[
 u=\frac{t}{t+2},\qquad x=3(t+2)^2X,\qquad
 y=\sqrt3(t+2)^3Y
\]
transforms the surface equation into $Y^2=X^3-2uX+u^4$. This isomorphism of generic elliptic curves extends
  to the relatively minimal elliptic surfaces. The quadratic twist cocycle is negation. On the self-product
  it is the diagonal Kummer involution, so it acts trivially on the quotient. Thus the quotient is isomorphic
  over $\mathbb F_5$ to the quotient in \cite[Example~9.1]{ref4}, which supplies the stated projective
  Calabi--Yau model. We choose $X_5$ to be this compatible model.

Let $K$ be the canonical crepant resolution of $\widehat Y/\langle\iota\rangle$, obtained by blowing up the
  image of the fixed curve. By Corollary~\ref{cor:betti-small} and Lemma~\ref{lem:quotient-cohomology},
\[
 b_3(K)=\dim H^3(\widehat Y_{\overline{\mathbb F}_5},\mathbb Q_\ell)^\iota
 +\dim H^1(C_{\overline{\mathbb F}_5},\mathbb Q_\ell)=0.
\]
The second term vanishes by Lemma~\ref{lem:fixed-curves}. Both $K$ and $X_5$ have trivial canonical bundle
  and are smooth proper birational threefolds. Corollary~\ref{cor:betti-birational} therefore gives
  $b_3(X_5)=0$.

For the canonical Kummer resolution $K$,
\[
 e(K)=\frac{e(\widehat Y)-e(C)}2+2e(C)
 =\frac{86+3\cdot14}{2}=64.
\]
The birational Calabi--Yau comparison gives the same Betti numbers for $X_5$. Since $b_1=b_3=b_5=0$, its
  second Betti number is $(64-2)/2=31$.

Suppose a smooth proper characteristic-zero lifting existed over a complete discrete valuation ring. An ample
  line bundle lifts because $H^2(X_5,\mathcal O_{X_5})=0$, so the lifting is projective after formal
  algebraization. Semicontinuity gives $h^1(\mathcal O)=h^2(\mathcal O)=0$ on the generic fiber; invariance
  of $\chi(\mathcal O)$ then gives $h^3(\mathcal O)=1$. Characteristic-zero Hodge theory forces $b_3\ge2$,
  whereas smooth proper base change gives $b_3=0$. This contradiction proves non-liftability to
  characteristic zero.
\end{proof}

\begin{proposition}
For the models above, Frobenius is zero on both $H^3(\widehat Y,\mathcal O_{\widehat Y})$ and
  $H^3(X_5,\mathcal O_{X_5})$. In particular $X_5$ is not Bloch--Kato $1$-ordinary, and hence is not
  $2$-ordinary.
\end{proposition}
\begin{proof}
First we justify the coherent comparison for $\rho:\widehat Y\to Y_+$. The only two-dimensional exceptional
  fiber is the smooth cubic surface $E$ at $\mathrm{IV}\times\mathrm{IV}$, with normal bundle $\mathcal
  O_E(-1)$. The cubic hypersurface sequence gives
\[
 H^j(E,\mathcal O_E(n))=0\qquad(j=1,2;\ n\ge0).
\]
Formal functions imply $R^2\rho_*\mathcal O=0$ there. The other fibers have dimension at most one. All
  positive higher direct images have finite support. Relative K\"unneth gives $H^1(Y_+,\mathcal
  O)=H^2(Y_+,\mathcal O)=0$. Since $R^2\rho_*\mathcal O=0$, Leray now gives $H^2(\widehat Y,\mathcal O)=0$.
  Crepancy and Serre duality give $H^1(\widehat Y,\mathcal O)=0$, and the low-degree Leray sequence forces
  $H^0(R^1\rho_*\mathcal O)=0$. Its support is finite, so $R^1\rho_*\mathcal O=0$. Therefore $R\rho_*\mathcal
  O_{\widehat Y}=\mathcal O_{Y_+}$, and K\"unneth yields
\[
 H^3(\mathcal O_{\widehat Y})\simeq H^1(\mathbb P^1,\mathcal O(-2)).
\]
For the cubic $f=4x^3+A(t)x+B(t)$ the Hasse coefficient in characteristic five is $[x^4]f^2=8A(t)$. In the
  present model it is
\[
 h(t)=4t(t+2)^3,\qquad h(t)^2=t^2(t+2)^6.
\]
On the \v{C}ech generator, Frobenius is $[t^{-1}]\mapsto[h(t)^2t^{-5}]$. Its only surviving coefficient is
\[
 [t^4]h(t)^2=\binom62 2^4=0\quad\text{in }\mathbb F_5.
\]
This proves the first assertion.

To pass to the Kummer model, let $\widehat Y_C=\operatorname{Bl}_C\widehat Y$. Its lifted involution has
  smooth quotient $K=\widehat Y_C/\langle\iota\rangle$, which is the canonical crepant resolution used above.
  The blowup $\beta:\widehat Y_C\to\widehat Y$ has $R\beta_*\mathcal O_{\widehat Y_C}=\mathcal O_{\widehat
  Y}$, and the quotient is tame. Hence
\[
 H^i(K,\mathcal O_K)\simeq H^i(\widehat Y_C,\mathcal O_{\widehat Y_C})^\iota
 \simeq H^i(\widehat Y,\mathcal O_{\widehat Y})^\iota
\]
compatibly with Frobenius. The action in degree three is $+1$ by Serre duality and the invariant volume form.
  In particular $K$ has $H^1(\mathcal O)=H^2(\mathcal O)=0$, and its top Frobenius is zero. The invariant
  volume form descends and extends without zeros by crepancy, so $\omega_K\simeq\mathcal O_K$.

Finally, global top forms on the smooth proper birational models $K$ and $X_5$ identify as rational forms.
  The Cartier operator on rational top forms is intrinsic to the common function field. By Serre duality it
  is dual to Frobenius on top coherent cohomology. Its vanishing on $K$ therefore implies its vanishing on
  $X_5$, proving the remaining assertion.
\end{proof}

The preceding propositions prove the Betti-number and Frobenius assertions of Theorem~\ref{thm:main}.

\section{The double-octic model and its lifting obstruction}

Unless otherwise specified, in this section \(k=\mathbb F_5\) and
\(W=W(k)=\mathbb Z_5\). In the general criterion of Lemma~\ref{lem:lam} we use
\(\kappa\) for an arbitrary perfect field of characteristic five.

\begin{theorem}\label{thm:double-octic}
The double octic below admits a smooth projective crepant resolution over
\(k\). Every such resolution \(X\) has no \(W_2(k)\)-lift and satisfies
\[
H^3_{\mathrm{cris}}(X/W)=0,\qquad
H^4_{\mathrm{cris}}(X/W)\simeq W^{31}.
\]
\end{theorem}

\subsection{The double-plane equation}

Put \(L=V+2U\). In \eqref{eq:surface}, write the affine base coordinate as \(t_0=V/U\)
and put \(x=t_0X/U\), \(y=t_0W_0/U^2\) on \(Ut_0\ne0\). After
cancelling \(t_0^2\), the resulting equation is
\[
W_0^2=A(U,V,X),\qquad A=4VX^3+3XL^3+3V^2L^2.
\]
Thus the double plane in \(\mathbb P(1,1,1,2)\) is birational over \(k\)
to the elliptic surface of the introduction, with negation acting by
\(W_0\mapsto-W_0\). On the chart \(L=1\), we use \(t=V/L\),
\(x_0=X/L\), and \(w=W_0/L^2\), so
\(w^2=4tx_0^3+3x_0+3t^2\). In these coordinates the fibers at
\(t=0,1,\infty\) have types \(\mathrm{III},I_5,\mathrm{IV}\),
respectively, by the change of base parameter.

\subsection{From the self fiber product to the double octic}

Let \(B=A(U,V,Y)\) and, on \(L\ne0\), put \(x_0=X/L\) and \(y_0=Y/L\).
On the generic fiber of the pencil the two double-plane coordinates satisfy
\(w_1^2=A/L^4\) and \(w_2^2=B/L^4\).
The simultaneous sign involution fixes their product. With \(t=V/L\),
its invariant function field is generated by \(t,x_0,y_0\) and
\(w_1w_2\). Consequently the quotient is birational to
\[
Z^2=F,\qquad F=AB,
\]
in \(\mathbb P(1,1,1,1,4)\), with \(Z=L^4w_1w_2\).
Thus the octic and the quotient \(X_+\) are birational over \(k\).
The covering involution is induced by negation in one factor.

For a blow-up of a smooth center of codimension \(c\) along which the
branch divisor has multiplicity \(m\), the new branch is the strict
transform plus the exceptional divisor if \(m\) is odd. The new covering
line bundle is the pullback of the old one minus \(\lfloor m/2\rfloor E\).
The double-cover canonical bundle formula gives discrepancy
\(c-1-\lfloor m/2\rfloor\) on the blown-up base. In particular the points
of multiplicities \(4,5\) and the double curves below are crepant centers
when used in an admissible resolution.\footnote{This is the blow-up rule
in \cite[\S 3]{ref6}. The multiplicities in the next table refer to the
original octic. They do not specify the incidence or order of all
infinitely near centers.}

\begin{center}
\begin{tabular}{clc}
 & ideal in \(k[U,V,X,Y]\) & \(m\)\\
\(P_1,P_2\) & \((U,V,X),(U,V,Y)\) & \(5\)\\
\(R_3,R_2\) & \((U,V-X,V-Y),(V-3U,X,Y)\) & \(4\)\\
\(L_{X,3},L_{Y,3}\) & \((U,V-X),(U,V-Y)\) & \(2\)\\
\(L_{X,2},L_{Y,2}\) & \((V-3U,X),(V-3U,Y)\) & \(2\)\\
\(L_\infty\) & \((U,V)\) & \(2\)
\end{tabular}
\end{center}

The two fivefold points consist of the quartic cone vertex of one factor
together with a smooth point of the other. The fourfold points occur
where the singular lines of the two factors at the same base value meet.
The double line \(L_\infty\) lies in both factors. The remaining
intersection is described by
\[
B-A=(Y-X)Q,\qquad Q=4V(X^2+XY+Y^2)+3L^3.
\]
The diagonal has ideal \((X-Y,A)\), and the off-diagonal ideal used in the
computation is \((A,Q):(U,V)^\infty\). Both occur with generic
multiplicity two in the branch divisor. The diagonal is birational to
the plane quartic \(A=0\).

\subsection{Construction and projectivity}

First blow up the two fivefold points and the fourfold points \(R_2,R_3\).
Blow up the further fourfold point over \(R_3\). Next blow up the two
successive exceptional curves over \(R_2\) and the two exceptional curves
over \(R_3\). Resolve the strict transforms of the two original singular
lines over each of \(R_2,R_3\) by one and two blow-ups per line,
respectively. Blow up three successive double curves over \(L_\infty\).
Finally separate the strict transforms of the diagonal and off-diagonal
curves and the four intersections with the branch components over the
fivefold points. On the resulting smooth base \(Y\), the branch divisor
has four ordinary double points, all over \(R_3\). Write
\(X_{\mathrm{nod}}\) for the double cover of \(Y\).

\medskip
\noindent\textit{Charts over the fourfold points.}
At \(R_3\), use \(s=t-1\) and complete the two branch factors in the
transverse coordinates. The common critical value has order five in
\(s\). Indeed, if \(c=(1+s)^{-1/2}\) is the critical section through
\(x_0=1\), the value is \(2c+3(1+s)^2\), and its order is five by
multiplying by the unit \(2c-3(1+s)^2\). Divide both factors by the same
leading transverse coefficient. The remaining transverse units have
square roots with residue one. The resulting branch factors are
\(a^2+s^5\epsilon(s)\) and \(b^2+s^5\epsilon(s)\), up to a common unit
whose product is a square. Here \(\epsilon(0)\ne0\). On the chart
\(a=s\alpha,b=s\beta\) the first point blow-up gives the factors
\(\alpha^2+s^3\epsilon(s)\) and \(\beta^2+s^3\epsilon(s)\). Blowing up
their common fourfold point leaves, on the \(s\)-chart, two smooth
factors of the form \(P\) and \(P+rz\), where \(r=v-u,z=v+u\) for
\(a=s^2u,b=s^2v\). Blowing up \((P,r)\) gives the branch \(p(p+z)\) on
\(P=rp\), and its remaining double curve is smooth. On the chart
\(s=\alpha\sigma,\beta=\alpha\gamma\) of the second point blow-up the
nonunit factor is \(\gamma^2+\alpha\sigma^3\epsilon(\alpha\sigma)\).
The exceptional curve blow-up \(\gamma=\sigma\eta\) leaves
\(\eta^2+\alpha\sigma\epsilon(\alpha\sigma)\), an ordinary double point
of the branch surface. At the other end, on the chart \(s=a\tau,b=a\beta\)
of the first point blow-up, the nonunit factor starts as
\(\beta^2+a^3\tau^5\epsilon(a\tau)\). After
\(\beta=a\gamma,\gamma=\tau^2\eta\) the corresponding factor is
\(\eta^2+a\tau\epsilon(a\tau)\). Interchanging \(a,b\) supplies the
other two nodes.

At \(R_2\), work on \(V=1\). The first point blow-up chart
\(Y=y,X=yv,L=yw\) gives factors
\[
3w^2+4yv^3+3y^2vw^3,\qquad 3w^2+4y+3y^2w^3.
\]
Blowing up \((y,w)\) and taking \(y=wq\) gives two smooth factors with
difference divisible by \(q\). The quotient of the difference restricts
to \(4(v^3-1)\) along the exceptional double curve. Its three roots are
simple. Near each root the branch has the form \(P(P+qr)\), and blowing
up \((P,q)\) separates the exceptional and residual double curves. On
the complementary chart \(w=yk\) the second factor is a unit near
\(y=k=v=0\) and the first is \(3yk^2+4v^3+3y^4vk^3\). Blowing up
\((k,v)\) gives a smooth branch there: in the charts \(v=kj\) and
\(k=vj\), the transformed equations have unit derivatives with respect
to \(y\) and \(v\), respectively, along the exceptional divisor.

\medskip
\noindent\textit{Charts along \(L_\infty\) and at its endpoints.}
Away from the fivefold points take \(X=1,Y=y\ne0\). The three successive
line blow-ups are seen in the charts \(V=L^jv,j=1,2,3\), by dividing each
factor by \(L^j\). On the last exceptional divisor the two factors
restrict to \(4v+3\) and \(4y^3v+3y\). They coincide only at \(y=\pm1\),
where their difference after eliminating \(v\) has a simple zero in \(y\).
The branch factors therefore meet transversely there. In the other
chart of each blow-up the transforms miss the exceptional direction
complementary to the displayed chart, since the coefficients of \(V\)
in both factors are units.

At \(X=1,L=V=Y=0\), use the point blow-up chart \(Y=y,L=y\ell,V=yv\).
The exceptional branch component is \(y=0\). After removing the
respective exceptional powers from \(A,B\), the two factors are
\[
4v+3y^2\ell^3+3y^3v^2\ell^2,\qquad 4v+3\ell^3+3v^2\ell^2.
\]
The line transforms are resolved in the charts \(v=\ell^jw\). After
the third blow-up their restrictions to \(\ell=0\) are \(4w+3y^2\)
and \(4w+3\). At \(y=0\) these give distinct roots and each branch meets
the exceptional branch plane transversely.

On the fivefold-point exceptional plane away from this direction, the
traces of the two factors are the line \(V=0\) and the quartic
\(4VY^3+3YL^3+3V^2L^2=0\) in the direction coordinates. They meet at
\([L:V:Y]=[1:0:0]\), where the line, the quartic, and the exceptional
branch divisor form a normal crossing. Blowing up the intersection of
the two quartic branch surfaces separates their intersections with the
exceptional plane there. At \(V\ne0\) only the quartic trace meets the
plane. Its two singular directions, on \(V=1\), are
\((L,Y)=(0,0)\) and \((1,1)\). The branch equation of this factor is
independent of the radial coordinate, so the line blow-ups resolving the
\(R_2\) and \(R_3\) cuspidal directions leave the transformed factor
smooth and transverse to the exceptional plane. Exchanging \(X\) and
\(Y\) gives the charts at the other fivefold point.

\begin{lemma}
The nodal double cover admits a smooth projective small resolution
\(X_{\mathrm{bl}}\).
\end{lemma}

\begin{proof}
The base \(Y\) is projective, since it is obtained by successive
blow-ups of \(\mathbb P^3_k\). The double cover \(X_{\mathrm{nod}}\to Y\)
is finite, so \(X_{\mathrm{nod}}\) is projective over \(k\).

Let \(E_a\) and \(E_b\) be the exceptional divisors of the blow-ups of
the two exceptional curves over \(R_3\). The cover splits over the
generic point of each divisor: in the charts of those blow-ups, the
restriction of the transformed branch equation to the exceptional
divisor is a square. Choose one of the resulting prime Weil divisors
on \(X_{\mathrm{nod}}\) over each of \(E_a\) and \(E_b\). Each chosen
divisor contains the two nodes over the corresponding exceptional
divisor. At either node the completed local equation is \(xy=uv\), and
the chosen Weil divisor has ideal \((x,u)\) after relabelling the
coordinates. The blow-up of \((x,u)\) is covered by the charts
\(u=xt,y=tv\) and \(x=us,v=sy\). Both charts are smooth, and the fibre
over the node is a single \(\mathbb P^1_k\). Away from the nodes the
chosen Weil divisors are Cartier, so the blow-ups are isomorphisms there.

Blowing up the two global ideal sheaves therefore gives a small crepant
resolution \(X_{\mathrm{bl}}\to X_{\mathrm{nod}}\). The target is
Gorenstein, so the small morphism is crepant. It is projective, as a
composition of blow-ups, and hence \(X_{\mathrm{bl}}\) is projective
over \(k\). The choice of a Weil component over each exceptional divisor
determines the rulings at both nodes on that divisor.
\end{proof}

\begin{proposition}\label{prop:nonliftable}
The threefold \(X_{\mathrm{bl}}\) does not lift to \(W_2(k)\).
\end{proposition}

\begin{proof}
Write a lift of the branch equation as \(\widehat F+5G\), where
\(\widehat F\) is the product of the indicated integer quartics and
\(G\in k[U,V,X,Y]_8\). At an original linear centre of multiplicity
\(m\), expand in lifted normal coordinates and retain the terms of
normal degree less than \(m\). Dividing the contribution from
\(\widehat F\) by \(5\) gives the inhomogeneous terms. The linear terms
in the centre motions are obtained by the chain rule. At subsequent
centres apply the same test to the branch transform, using a common
motion whenever a centre enters more than one test. Over \(R_3\) the
conditions include the fourfold point jet and the node jets
\([s^4a^2]\) and \([s^4b^2]\). Over \(R_2\), after the first exceptional
curve has been moved, write the local equation as \(P(P+qr)+5g\).
Regularity of the next centre requires
\[
g|_{P=q=0}=0,\qquad r\mid g_q,\qquad r^2\mid rg_P-2g_q,
\]
with the derivatives evaluated along the centre. Indeed, if its normal
motion is \((p,d)\), the linear terms require \(g_q+rp=0\) and
\(g_P+2p+rd=0\). The same motion of \(R_2\) is used for the successive
jet conditions.

\medskip
\noindent\textit{The node rows over \(R_3\).}
Use \(s=t-1,a=x_0-c(s),b=y_0-c(s)\), where \(c(s)\) is the common
critical section. Consider the node on the chart with exceptional
parameter \(a\). The successive substitutions there have composite
\[
s=a\tau,\qquad b=a^2\tau^2\eta.
\]
The branch equation is divided by \(a^6\tau^4\): \(a^4\) comes from the
fourfold point, \(a^2\) from the exceptional double curve, and
\(\tau^4\) from the two blow-ups of the line. The constant term in the
resulting node chart is therefore read from \(s^4a^2\). Indeed, for
\(s^ia^jb^k\) to contribute at \(a=\tau=\eta=0\) one needs \(k=0,i=4\),
and \(i+j=6\). The symmetric chart gives \(s^4b^2\).

The centre motions can be checked in the normalized special-fibre
equations, since their contributions are multiplied by \(5\). For the
first point blow-up a translation in \(b\) differentiates a factor with
vanishing \(b\)-derivative at \(b=0\). Translations in \(s,a\) likewise
contribute nothing to \(s^4a^2\): at \(b=0\) the equation has the form
\(s^5\epsilon(s)(a^2+s^5\epsilon(s))\), up to a unit, and
\(\partial_s(s^5)=0\) in \(k\). For the exceptional curve
\((a,\beta)=(0,0)\), use the special-fibre factor
\(\beta^2+a^3\tau^5\epsilon(a\tau)\). Translation in \(\beta\)
contributes a term linear in \(\eta\). Translation in \(a\)
contributes, after the remaining divisions, a multiple of \(\tau\)
or \(a\). For the two line blow-ups the corresponding equations, up to
units, are \(\gamma^2+a\tau^5\epsilon(a\tau)\) and
\(\delta^2+a\tau^3\epsilon(a\tau)\). Translation in the squared
coordinate contributes a term vanishing at \(\eta=0\). Translation
in \(\tau\) contributes a term divisible by \(a\) after the remaining
divisions. Unit changes do not affect the vanishing of the node
parameter. Finally a translation of the node itself contributes zero
to the constant term since the special-fibre gradient vanishes there.
Thus these two rows do not require a choice of lifts of the centres.

\medskip
\noindent\textit{The endpoint conditions over \(L_\infty\).}
Work over the open line with \(X=1,Y=y\ne0\). The branch equation has a
nondegenerate quadratic term in \(V\), so we may complete the square
formally along its critical section. The critical section modulo
\(L^6\) is \(V=-3(1+y^{-2})L^3/8\), obtained from the linear and
quadratic terms in \(V\). The critical value on the special fibre is
\(L^6u(L,y)\), with
\[
u(0,y)=-\frac94 y^3(1-y^{-2})^2.
\]
Dividing the branch equation by the unit \(y^3\) gives the normalized
coefficient \(-9(1-y^{-2})^2/4\). For the unnormalized equation the
jet conditions are imposed on \(G\) as written; multiplication by
this unit does not change the endpoint divisibilities or evaluations.
A first-order motion in \(V\) contributes nothing to the critical
value, since its \(V\)-derivative vanishes there.

Each of the three double-curve blow-ups divides this value by the
square of a lifted local equation for the exceptional parameter. On
the open line, replacing \(L\) by \(L+5d(y)\) at any stage changes the
coefficient one degree below the leading order by a regular multiple
of \(u(0,y)\). More explicitly, when the remaining critical value
starts in degree \(2j\), the translation contributes
\(5(2j)d(y)u(0,y)L^{2j-1}\). Pulling this contribution back through
the preceding divisions gives a term of degree \(5\) in the original
critical value. Thus all corrections below degree \(5\) must vanish,
and the degree-five coefficient must lie in \((u(0,y))\) near each
crossing. The moving centres are regular over both \(y=1\) and
\(y=-1\), so this is a double-zero condition at each point.

For the degree-six coefficient, the two reduced branch roots after
the third division coincide at \(L=0,y=\pm1\). The special-fibre
critical value after division is, up to a unit, the square of their
difference. Its first derivatives in \(L\) and \(y\) vanish at these
two points. Consequently changing the lifted coordinates or moving
the intersection curve does not change its first-order value there.
A lifted double curve through the crossing forces that value to
vanish, giving the two degree-six evaluation rows.

Let \(M_{\mathrm{lift}}\) be the coefficient matrix of the affine
lifting system and let \(b\) be its right-hand side. The rank
computation gives
\[
\operatorname{rank}M_{\mathrm{lift}}=149,\qquad
\operatorname{rank}(M_{\mathrm{lift}}\mid b)=150.
\]
The rank discrepancy excludes a lift of the marked branch data. We
now show that a lift of \(X_{\mathrm{bl}}\) would supply such data.
For the small contraction \(X_{\mathrm{bl}}\to X_{\mathrm{nod}}\),
the pushforward of the structure sheaf is the structure sheaf and the
first higher direct image vanishes. Thus any lift descends to a lift
of \(X_{\mathrm{nod}}\) with a singular section through every node,
by \cite[Theorems 3.1 and 4.1 and Proposition 4.2]{ref7}.

The classes of node-preserving lifts, if nonempty, form an affine
subspace of the torsor of liftings: they are cut out by the local node
parameters. This subspace is stable under the covering involution.
Taking the midpoint of a class and its translate gives an invariant
class. A lift of the involution on a representative can be made
involutive by composing it with minus half the infinitesimal
automorphism given by its square. The quotient \(\mathcal Y\) is a
flat lift of \(Y\), since taking invariants is a direct summand and
commutes with reduction. The finite double-cover algebra is the sum
of its invariant structure sheaf and an invertible anti-invariant
summand. Indeed, the latter is flat over the Artinian base and reduces
to a line bundle. Its multiplication therefore gives a lifted branch
section.

For clarity, write \(f:Y'\to Y_0\) for one blow-up in the construction
of \(Y\), with smooth centre \(C\) and exceptional divisor \(E\).
Given a lift of \(Y'\), the vanishings
\(H^0(E,\mathcal O_E(E))=H^1(E,\mathcal O_E(E))=0\) lift \(E\)
uniquely as a Cartier divisor. They follow from
\(\mathcal O_E(E)=\mathcal O_{\mathbb P(N_{C/Y_0})}(-1)\).
The structure-sheaf blow-down lifts \(f\). Pushing forward
\[
0\longrightarrow\mathcal O_{\mathcal Y'}(-\mathcal E)
\longrightarrow\mathcal O_{\mathcal Y'}
\longrightarrow\mathcal O_{\mathcal E}\longrightarrow0
\]
is exact and commutes with reduction, since the first higher direct
images of all three special-fibre sheaves vanish. For the first
sheaf this follows from the same sequence on the special fibre and
the surjection \(\mathcal O_{Y_0}\to\mathcal O_C\); for the other two
it follows from the blow-up and projective bundle formulas. The
quotient \(\widetilde f_*\mathcal O_{\mathcal E}\) is flat over
\(W_2(k)\) by the square-zero base-change sequence. Consequently the
kernel is the ideal of a flat lifted centre \(\mathcal C\). The
pullback of this ideal generates
\(\mathcal O_{\mathcal Y'}(-\mathcal E)\) by Nakayama's lemma. The
universal property gives a morphism to
\(\operatorname{Bl}_{\mathcal C}\mathcal Y_0\), which is an
isomorphism by flatness and reduction. Here the lifted base and
centre are smooth over \(W_2(k)\), since they are flat with smooth
special fibres.

The components of the final branch divisor are disjoint, so they
lift separately as open-and-closed parts of the lifted divisor. At
each blow-down track their lifted equations individually. A branch
component equal to \(E\) on the special fibre is the unique lift
\(\mathcal E\) and is dropped in the blow-down. For each other
component, let \(D_0\) be its special-fibre image and let \(m_0\) be
its multiplicity at \(C\). Its strict transform has line bundle
\(f^*\mathcal O(D_0)(-m_0E)\). All the intermediate bases have
\(H^1(\mathcal O)=H^2(\mathcal O)=0\), so this identity lifts, with
\(\mathcal O(D_0)\) replaced by its unique lift. Pushing the section
down and using
\(\widetilde f_*\mathcal O(-m_0\mathcal E)=I_{\mathcal C}^{m_0}\)
gives a branch equation with the required multiplicity. Iterating
this procedure recovers the marked branch transforms down to
\(\mathbb P^3_{W_2(k)}\). The resulting equation corrections and
centre motions must satisfy the affine system above, which is
inconsistent.
\end{proof}

\begin{lemma}\label{lem:rigidity-double}
The threefold \(X_{\mathrm{bl}}\) satisfies
\[
H^1(X_{\mathrm{bl}},\Omega^2_{X_{\mathrm{bl}}/k})
=H^2(X_{\mathrm{bl}},\Omega^1_{X_{\mathrm{bl}}/k})=0.
\]
\end{lemma}

\begin{proof}
Let \(\pi:X_{\mathrm{nod}}\to Y\) be the double cover, with covering
line bundle \(\mathcal L\simeq\omega_Y^{-1}\).
The double-cover tangent-sheaf calculation of \cite[\S 2]{ref6} also
applies to this nodal cover, since \(2\) is invertible in \(k\).
Indeed, locally for \(w^2=f(z_1,z_2,z_3)\), an anti-invariant
derivation has the form
\[
D(z_i)=w b_i,\qquad
D(w)=\frac12\sum_i b_i\,\frac{\partial f}{\partial z_i}.
\]
Thus
\[
(\pi_*T_{X_{\mathrm{nod}}})^-
\simeq T_Y\otimes\mathcal L^{-1}
\simeq\Omega_Y^2.
\]
The anti-invariant part of the locally trivial nodal tangent space is
therefore \(H^1(Y,\Omega_Y^2)=0\), since the curve centres are
geometrically rational. In particular the diagonal is the normalization
of a quartic with cusp delta invariants \(1,2\); the off-diagonal curve
parametrizes ordered pairs of distinct roots of the cubic, with inertia
orders \(2,2,3\). Its degree-six cover has genus zero by
Riemann--Hurwitz.

For the invariant part, let \(M_{\mathrm{eq}}\) denote the matrix of
the homogeneous system with the restrictions
\(G|_{C_\Delta}=G|_{C_{\mathrm{off}}}=0\).
The computation gives \(\operatorname{rank}M_{\mathrm{eq}}=150\).
The restrictions follow by testing at the generic points of the two
double curves, where both \(F\) and \(dF\) vanish.
The sixteen polynomials \(Z_i\partial F/\partial Z_j\), with
\(Z=(L,V,X,Y)\), span a space of dimension \(15\) over \(k\).
This is a separate rank check on their degree-eight coefficient vectors.
They lie in the homogeneous kernel, since a linear coordinate change
moves the entire marked configuration. They therefore exhaust the
kernel. After a coordinate change we may set \(G=0\).

We check the remaining motions in the order of the blow-ups. At a
fivefold point, on \(X=1\), the degree-five initial form is
\(4V(4VY^3+3YL^3+3V^2L^2)\). Its three first partial derivatives in
\(L,V,Y\) are independent. A translation preserving multiplicity
five must therefore be zero. The other point is symmetric.

At \(R_3\), the initial form \(a^2b^2\), up to a unit, fixes the two
transverse coordinates of the first point. The weight vector field of
weights \((0,1,2,2)\) on \((L,V,X,Y)\) preserves the branch divisor
and moves \(t\) with nonzero velocity at \(1\). It thus fixes the
remaining freedom in the location of \(R_3\) after a coordinate
change. At the next fourfold point the same initial form again
fixes the two transverse coordinates. If the point moves by
\(s=\varepsilon d\), the chart
\(s=\varepsilon d+\alpha\sigma,\beta=\alpha\gamma\) has nonunit factor
\[
\gamma^2+\alpha\sigma^3\epsilon(\alpha\sigma)
+3\varepsilon d\,\sigma^2\epsilon(0)
\quad\text{up to terms giving zero at the node.}
\]
Dividing by \(\sigma^2\) on the ensuing exceptional-curve chart shows
that the node parameter changes by \(3\varepsilon d\,\epsilon(0)\).
Moving the curve does not cancel this constant: motion in \(\gamma\)
gives a term vanishing at \(\eta=0\), and motion in \(\sigma\) gives
a term divisible by \(\alpha\). Hence \(d=0\). Each exceptional
curve over \(R_3\) projects to a line in the first exceptional plane.
Its projected motion out of that plane lies in
\(H^0(\mathbb P^1,\mathcal O(-1))=0\). At a generic point the
quadratic term in the remaining normal coordinate fixes the motion
in the plane. The regular normal section then vanishes on the whole
curve.

At \(R_2\), the initial form \(9L^4\) fixes the \(L\)-coordinate of
the point. Write its other translations as \(\varepsilon p\) and
\(\varepsilon q\). In the chart \(X=a,Y=av,L=aw\), the lowest terms
of the two factors are \(3w^2+4a+12\varepsilon p\) and
\(3w^2+4av^3+12\varepsilon qv^2\). If the first double curve moves in
the \(a\)-direction by \(\varepsilon d\), its linear \(a\)
coefficient imposes
\[
32v^3d+48(pv^3+qv^2)=0.
\]
After the curve blow-up \(a=wz\), the constant correction on the
second exceptional double curve is read from the coefficient of
\(w^2\) before division. Substituting \(d\) gives
\(18(1-v^3)(p-q/v)\). A double curve requires this to vanish at its
generic point, so \(p=q=0\).

The first exceptional curve over \(R_2\) can now move only within
the exceptional plane. Its normal motion there is a section of
\(\mathcal O_{\mathbb P^1}(1)\). Near an end use the chart
\(w=yk+\varepsilon c(v)\) of the first curve blow-up. The factor
meeting the original singular line becomes
\(3yk^2+4v^3+6\varepsilon c(v)k+\cdots\). Preserving the double
centre \(k=v=0\) requires the motion of \(k\) at its generic point
to be \(-\varepsilon c(0)/y\). Regularity at \(y=0\) forces
\(c(0)=0\). The opposite end gives the second zero by symmetry.
Thus the first curve is fixed, and the generically nondegenerate
normal quadratic form fixes the second exceptional curve.

The centres over the original singular lines and over \(L_\infty\)
project isomorphically to those lines. Their projected normal motion
is a section of \(\mathcal O_{\mathbb P^1}(1)^{\oplus2}\). It
vanishes at both ends: the lifted centre meets transversely the
exceptional divisor over each fixed endpoint, and this intersection
maps to that endpoint. For the first centre over an original line,
projection is an isomorphism at its generic point. For the second
centre over a line at \(t=1\), the projected motion fixes the base
direction at the generic point and the unit quadratic coefficient
fixes the other normal coordinate. Likewise, on the successive
\(L_\infty\) charts \(v=\ell^jw\) the projection retains the
\(\ell\)-direction and the unit quadratic coefficient fixes the
\(w\)-direction. This gives zero motion on each curve in turn.
Finally, at a generic point of each remaining intersection curve
the local branch equation is a product of two transverse parameters.
Its linear coefficients after translation force both normal motions
to vanish. In each case a regular normal section that is zero
generically is zero everywhere.

This leaves no invariant tangent direction. Restriction of
\(H^1(T_{X_{\mathrm{bl}}})\) off the exceptional curves is injective by
depth, so the descent to the locally trivial nodal tangent space also
detects every first-order deformation of the smooth model; compare
\cite[Theorem~4.1 and Proposition~3.2]{ref7}.
\end{proof}

The following criterion is used by Lam
\cite[Corollary~1.6 and Lemmas~5.1--5.2]{ref16}.
We recall the argument in characteristic five using
\cite[Theorem~1.3 and Corollary~4.4]{ref1}.

\begin{lemma}\label{lem:lam}
Let \(T\) be a smooth projective Calabi--Yau threefold over a
perfect field \(\kappa\) of characteristic \(5\).
Suppose
\(H^1(T,\Omega^2_{T/\kappa})=H^2(T,\Omega^1_{T/\kappa})=0\)
and \(T\) has no \(W_2(\kappa)\)-lift.
Then \(H^3_{\mathrm{dR}}(T/\kappa)=0\).
\end{lemma}

\begin{proof}
Put \(T^{(1)}=T\times_{\kappa,F_\kappa}\kappa\).
We use the conjugate spectral sequence in the convention
\[
\begin{aligned}
{}^{\mathrm{conj}}E_2^{i,j}
=H^i(T^{(1)},\Omega^j_{T^{(1)}/\kappa})
&\ \Longrightarrow\ H^{i+j}_{\mathrm{dR}}(T/\kappa),\\
d_2^{i,j}:E_2^{i,j}&\longrightarrow E_2^{i+2,j-1}.
\end{aligned}
\]
Since \(\kappa\) is perfect, \(T^{(1)}\) also has no
\(W_2(\kappa)\)-lift. Let
\(\theta\in H^2(T^{(1)},T_{T^{(1)}/\kappa})\)
be its nonzero lifting obstruction. By
\cite[Theorem~1.3 and Corollary~4.4]{ref1}, \(d_2\) is, up to sign,
cup product followed by contraction with \(\theta\).
All Hodge groups in the following two maps are on \(T^{(1)}\):
\[
H^1(\Omega^1)\longrightarrow H^3(\mathcal O),\qquad
H^0(\Omega^3)\longrightarrow H^2(\Omega^2).
\]
The first map is surjective by the perfect Serre pairing with
\(H^2(T_{T^{(1)}/\kappa})\). Contraction with a volume form identifies
\(T_{T^{(1)}/\kappa}\) with \(\Omega^2_{T^{(1)}/\kappa}\), so the
second map is injective. By hypothesis the two mixed terms of total
degree three vanish, and these maps kill the remaining two terms.
Hence \(H^3_{\mathrm{dR}}(T/\kappa)=0\).
\end{proof}

\begin{proposition}\label{prop:crystalline-middle}
For \(W=W(k)\) we have \(H^3_{\mathrm{cris}}(X_{\mathrm{bl}}/W)=0\)
and \(H^4_{\mathrm{cris}}(X_{\mathrm{bl}}/W)\simeq W^{31}\).
\end{proposition}

\begin{proof}
By the canonical bundle formula for the double cover and crepancy,
\(X_{\mathrm{bl}}\) has trivial canonical bundle. The double-cover
decomposition \(\pi_*\mathcal O=\mathcal O_Y\oplus\omega_Y\), the
vanishing \(H^i(Y,\mathcal O_Y)=0\) for \(i>0\), and
\(R^if_*\mathcal O_{X_{\mathrm{bl}}}=0\) for \(i>0\) for the small
nodal contraction \(f\), give the Calabi--Yau vanishings. Since
\(X_{\mathrm{bl}}\) and \(X_5\) are birational over \(k\), Corollary~\ref{cor:betti-birational}
and Proposition~\ref{prop:betti-kummer} give \(b_4(X_{\mathrm{bl}})=31\).
Proposition~\ref{prop:nonliftable} and Lemma~\ref{lem:rigidity-double} give
\(H^3_{\mathrm{dR}}(X_{\mathrm{bl}}/k)=0\) by Lemma~\ref{lem:lam}.
The crystalline base-change sequence
\[
0\longrightarrow H^3_{\mathrm{cris}}(X_{\mathrm{bl}}/W)/5
\longrightarrow H^3_{\mathrm{dR}}(X_{\mathrm{bl}}/k)
\longrightarrow H^4_{\mathrm{cris}}(X_{\mathrm{bl}}/W)[5]
\longrightarrow0
\]
gives the vanishing of \(H^3_{\mathrm{cris}}\) by Nakayama's lemma and
the freeness of \(H^4_{\mathrm{cris}}\).
Its rank is \(b_4\) by the crystalline--\'etale comparison of
Frobenius polynomials \cite{ref15}.
\end{proof}

\begin{lemma}\label{lem:crystalline-birational}
Let \(X_1\) and \(X_2\) be birational smooth projective threefolds over
a perfect field of characteristic five, with trivial canonical bundles.
The closure of the graph of a birational map and its transpose induce
inverse integral crystalline correspondences.
\end{lemma}

\begin{proof}
By Proposition~\ref{prop:codim-one}, there is a common open \(U\) whose complements
have codimension at least two. Let \(\Gamma\) be the closure of the
graph of the isomorphism. In the refined intersection defining
\({}^t\Gamma\circ\Gamma\), the part over \(U\) is the diagonal with
multiplicity one. The other contributions push forward into
\((X_1\setminus U)\times(X_1\setminus U)\), of dimension at most two.
They have zero pushforward as cycles of dimension three.
Thus the composition is \(\Delta_{X_1}\) in the integral Chow group,
and the reverse composition is \(\Delta_{X_2}\).
Integral crystalline cycle classes and the projection formula give
inverse crystalline correspondences, including on torsion \cite{ref11}.
In particular the two crystalline groups just computed are the same
for every such resolution.
\end{proof}

\begin{proof}[Proof of Theorem~\ref{thm:double-octic}]
Proposition~\ref{prop:crystalline-middle} and Lemma~\ref{lem:crystalline-birational} give the asserted
  crystalline groups.
Finally, for any such resolution \(X\), the crystalline reduction
sequence gives \(H^3_{\mathrm{dR}}(X/k)=0\).
A \(W_2(k)\)-lift would give \(E_1\)-degeneration by
Deligne--Illusie \cite{ref9} because \(3<5\), contradicting
\(H^0(X,\Omega_X^3)=k\). Thus non-liftability also holds independently
of the chosen resolution. It also applies to \(X_5\), using the
birational identification over \(k\) above and the triviality of both
canonical bundles.
\end{proof}

\subsection{Reproducing the rank computations}

The supplementary files \texttt{verify\_double\_octic\_ranks.m}
and \texttt{verify\_double\_octic\_ranks.sage} implement the affine and
homogeneous systems used above.
Here \(M_{\mathrm{lift}}\) denotes the affine coefficient matrix
from Proposition~\ref{prop:nonliftable}, with right-hand side \(b\), whereas
\(M_{\mathrm{eq}}\) denotes the homogeneous matrix in Lemma~\ref{lem:rigidity-double},
with the restrictions \(G|_{C_\Delta}=G|_{C_{\mathrm{off}}}=0\).
The rank statements are
\[
\operatorname{rank}M_{\mathrm{lift}}=149,\qquad
\operatorname{rank}(M_{\mathrm{lift}}\mid b)=150,\qquad
\operatorname{rank}M_{\mathrm{eq}}=150.
\]
The companion notes describe the ordering of the columns and the
blocks of rows. The implementation uses the supplied generators of
the saturated off-diagonal ideal.

\section{Crystalline cohomology and deformation consequences}

In this section \(k=\overline{\mathbb F}_5\) and \(W=W(k)\).
We write \(X_5\) also for its base change
\(X_5\times_{\mathbb F_5}k\).
Deformations are marked by an identification of their special fiber
with the given variety.

\begin{lemma}\label{lem:global-forms}
The projective Calabi--Yau threefold \(X_5\) satisfies
\[
 H^0(X_5,\Omega^1_{X_5/k})=H^0(X_5,\Omega^2_{X_5/k})=0.
\]
In particular, \(H^0(X_5,T_{X_5/k})=0\).
\end{lemma}

\begin{proof}
Put \(S=S_+\times_{\mathbb F_5}k\) and
\(Z=S\times_{\mathbb P^1_k}S\), with projections \(p_i\)
and structure map \(g\).
Write \(\pi:S\to\mathbb P^1_k\) for the elliptic fibration and
\(O\) for its zero section.
Let \(\rho:\widehat Y\to Z\) be the base change of the equivariant
crepant resolution constructed in Lemma~\ref{lem:equivariant}.
The singularities of \(Z\) are isolated hypersurface singularities.

The conormal presentation shows that \(\Omega_Z^1\) has depth at least two at these points.
It is torsion-free, since a torsion subsheaf would be supported at the isolated singularities and would have
  depth zero.
Consequently \(\Omega_Z^1\) is reflexive.
Restriction to the smooth locus and pullback to \(\widehat Y\) therefore give
\[
 H^0(\widehat Y,\Omega^1_{\widehat Y})=H^0(Z,\Omega_Z^1).
\]
The relative conormal sequence is
\[
 0\longrightarrow g^*\Omega_{\mathbb P^1}^1
 \longrightarrow p_1^*\Omega_S^1\oplus p_2^*\Omega_S^1
 \longrightarrow\Omega_Z^1\longrightarrow0.
\]
Since \(S\) is rational and \(p_{i*}\mathcal O_Z=\mathcal O_S\), the middle term has no global sections.
Relative K\"unneth gives \(H^1(Z,g^*\Omega^1_{\mathbb P^1})=k\).
Its generator maps, in the first summand, to \(p_1^*c_1(F)\), where \(F\) is a fibre of \(S\to\mathbb P^1\).
This class is nonzero: \(p_1\) has the section \(s\mapsto(s,O(\pi(s)))\), and restriction of \(c_1(F)\) to
  the zero section of \(S\) has degree one.
Thus \(H^0(\widehat Y,\Omega^1)=0\).

Next \(H^0(S,T_S)=0\).
Indeed, a vector field on \(S\) induces a vector field on the base through the anticanonical pencil.
At a critical point of the elliptic fibration the differential to the base is zero, so the induced vector
  field vanishes at each of the three distinct critical values.
Here we use critical points on the smooth minimal surface, not multiplicities in its discriminant; at an
  \(I_5\) node the local equation of the map is \(t=uv\), and \(dt\) vanishes at the node.
A vector field on \(\mathbb P^1\) with three zeros is zero.
The original vector field is therefore vertical.
The critical locus on \(S\) is finite, and the kernel of \(T_S\to\pi^*T_{\mathbb P^1}\) is the line bundle
\[
 T_{S/\mathbb P^1}\simeq\omega_S^{-1}\otimes\pi^*\omega_{\mathbb P^1}
 \simeq\mathcal O_S(-F),
\]
which has no sections.
The morphism \(\widehat Y\to S\times S\) has injective differential at its generic point.
Moreover \(R\rho_*\mathcal O_{\widehat Y}=\mathcal O_Z\), as proved above, and hence \((p_i\rho)_*\mathcal
  O_{\widehat Y}=\mathcal O_S\).
It follows that \(H^0(\widehat Y,T_{\widehat Y})=0\).
Contraction with the volume form gives \(H^0(\widehat Y,\Omega^2)=0\).

Let \(\widehat Y_C=\operatorname{Bl}_C\widehat Y\) and \(K=\widehat Y_C/\langle\iota\rangle\).
Blowing up a smooth centre of codimension two does not change global regular differential forms: restrict a
  form to the complement of the centre and use Hartogs extension on the smooth target.
Pullback under the generically separable double cover \(\widehat Y_C\to K\) is injective on forms.
Thus \(K\) has no global forms of degrees one and two.
Finally \(K\) and \(X_5\) are isomorphic in codimension one by Proposition~\ref{prop:codim-one}; Hartogs
  extension identifies their global differential forms.
This proves the assertions.
\end{proof}

\begin{theorem}\label{thm:rigidity}
The threefold \(X_5\) is infinitesimally rigid and has no lift to \(W_2(k)\).
\end{theorem}

\begin{proof}
The graph correspondence of Lemma~\ref{lem:crystalline-birational} identifies the middle
crystalline groups of \(X_5\) with those of \(X_{\mathrm{bl}}\),
computed in Proposition~\ref{prop:crystalline-middle}.
After extension of scalars to \(W\), crystalline base change gives
\(H^3_{\mathrm{dR}}(X_5/k)=0\).
Non-liftability follows by the same Deligne--Illusie argument as in
the proof of Theorem~\ref{thm:double-octic}.

For the conjugate spectral sequence we use the indexing
\[
 E_2^{p,q}
 =H^p\bigl(X_5^{(1)},\Omega^q_{X_5^{(1)}/k}\bigr)
 \ \Longrightarrow\ H^{p+q}_{\mathrm{dR}}(X_5/k).
\]
The term \(E_2^{1,2}\) has no incoming differential.
Its only possible outgoing differential is \(d_2\) to
\(H^3(X_5^{(1)},\Omega^1_{X_5^{(1)}/k})\).
By Serre duality, this target is dual to
\(H^0(X_5^{(1)},\Omega^2_{X_5^{(1)}/k})\), which vanishes
by Lemma~\ref{lem:global-forms}.
Thus \(E_2^{1,2}=E_\infty^{1,2}=0\).
Since \(k\) is perfect, \(H^1(X_5,\Omega^2_{X_5/k})=0\);
contraction with a volume form gives \(H^1(X_5,T_{X_5/k})=0\).
\end{proof}

\begin{corollary}
The Artin--Mazur formal group of \(X_5\) has infinite height.
In particular, \(X_5\) is not 1-ordinary or Bloch--Kato 2-ordinary.
\end{corollary}

\begin{proof}
By Katz--Messing~\cite{ref15} and Theorem~\ref{thm:main},
\(H^3_{\mathrm{cris}}(X_5/W)\otimes_W W[1/5]=0\).
The slope-\([0,1)\) comparison identifies
\(H^3(X_5,W\mathcal O_{X_5})\otimes_W W[1/5]\)
with the corresponding part of this isocrystal.
If the Artin--Mazur height were finite, say \(h\), this space would
have dimension \(h>0\); see~\cite[\S3]{ref25}.
Thus the height is infinite.
This also agrees with the explicit Hasse--Witt calculation above.
\end{proof}

\begin{theorem}
The crystalline cohomology of \(X_5\) is torsion-free.
As \(W\)-modules, with Frobenius and Tate twists omitted, it is
\[
 \bigl(H^i_{\mathrm{cris}}(X_5/W)\bigr)_{i=0}^6
 \simeq\bigl(W,0,W^{31},0,W^{31},0,W\bigr).
\]
The deformation tangent space \(H^1(X_5,T_{X_5})\) is zero, whereas \(H^2(X_5,T_{X_5})\) has dimension 32.
The universal formal deformation ring over \(W\) is \(k\).
\end{theorem}

\begin{proof}
The vanishing of the tangent space and the failure of a
\(W_2\)-lifting are proved in Theorem~\ref{thm:rigidity}.
Lemma~\ref{lem:global-forms} and \(H^1(\mathcal O)=0\) give
\(H^1_{\mathrm{dR}}(X_5/k)=0\).
Derived crystalline base change yields
\[
 0\longrightarrow H^1_{\mathrm{cris}}/5
 \longrightarrow H^1_{\mathrm{dR}}
 \longrightarrow H^2_{\mathrm{cris}}[5]
 \longrightarrow 0.
\]
Consequently \(H^1_{\mathrm{cris}}=0\) and
\(H^2_{\mathrm{cris}}\) is torsion-free.
Integral crystalline duality pairs torsion in degrees \(i\) and \(7-i\)
by \(W[1/5]/W\)-duality~\cite[Theorem~6.4 and Remark~6.6]{ref24}.
The middle groups are given by Proposition~\ref{prop:crystalline-middle} and
  Lemma~\ref{lem:crystalline-birational},
with scalar extension to \(W\).
The degree-zero and top groups are \(W\).
The remaining free ranks follow from Katz--Messing~\cite{ref15} and
Proposition~\ref{prop:betti-kummer}.

Riemann--Roch gives \(\chi(\Omega^1)=-e(X_5)/2=-32\).
The global-form vanishings, rigidity and Serre duality therefore give
\(h^1(\Omega^1)=h^2(T)=32\).
Finally, over a small extension with kernel \(I\), the kernel of
restriction on automorphisms of a deformation is
\(H^0(X_5,T_{X_5})\otimes_k I\).
Its vanishing implies by induction that every automorphism preserving
the marking is the identity.
Schlessinger's criterion therefore gives prorepresentability of the
marked deformation functor over \(W\)~\cite[Theorem~2.11 and Proposition~3.10]{ref23}.
Its complete local \(W\)-algebra \(R\) has no relative tangent
directions, so \(W\to R\) is surjective by the complete local
Nakayama lemma.
Thus either \(R=W\) or \(R=W/(5^n)\) for some \(n\geq1\).
The absence of a \(W_2\)-point excludes the first possibility and
forces \(n=1\) in the second. Hence \(R=k\).
\end{proof}

\begin{samepage}
\begin{corollary}
Write \(h^{p,q}=\dim_k H^q(X_5,\Omega^p_{X_5/k})\).
The only nonzero Hodge numbers of \(X_5\) are
\[
 h^{0,0}=h^{3,0}=h^{0,3}=h^{3,3}=1,
 \qquad h^{1,1}=h^{2,2}=32.
\]
With the indexing \(E_1^{p,q}=H^q(X_5,\Omega^p_{X_5/k})\),
its Hodge--de Rham spectral sequence satisfies \(E_1=E_2\) and \(E_3=E_\infty\).
Its only nonzero differentials have rank one and are
\[
 \begin{aligned}
 d_2^{1,1}&:H^1(\Omega^1)\longrightarrow H^0(\Omega^3),\\
 d_2^{0,3}&:H^3(\mathcal O)\longrightarrow H^2(\Omega^2).
 \end{aligned}
\]
\end{corollary}
\end{samepage}

\begin{proof}
The Hodge numbers follow from the preceding theorem, the Calabi--Yau conditions, Lemma~\ref{lem:global-forms}
  and Serre duality.
They force every \(d_1\) to vanish and leave only the displayed \(d_2\).
Since \(H^3_{\mathrm{dR}}=0\), both have rank one.
No higher differential can occur.
\end{proof}

\begin{corollary}
The threefold \(X_5\) is non-Hodge--Witt and its Hodge--Witt number \(h_W^{1,2}(X_5)\) is \(-1\).
\end{corollary}

\begin{proof}
We use Ekedahl's Hodge--Witt numbers~\cite{ref10}.
The Betti numbers give \(\deg c_3(X_5)=64\).
Joshi's formula~\cite[Proposition~7.2.1]{ref14} gives
\[
 h_W^{1,2}(X_5)
 =b_2(X_5)-\frac12\deg c_3(X_5)
 =31-32=-1.
\]
The non-Hodge--Witt assertion follows from~\cite[Theorem~7.3.1]{ref14}.
\end{proof}

\begin{remark}
The example combines torsion-free crystalline cohomology with a
negative Hodge--Witt number and failure of Hodge--de Rham degeneration.
The positions of the two nonzero differentials agree with those
computed for Hirokado's characteristic-three example~\cite[Remarks~1.7 and~5.9]{ref5}.
The equality \(b_3(X_5)=0\) also provides an instance of the vanishing
asked about in~\cite[Question~1.3]{ref16} for Calabi--Yau threefolds that do not
lift to characteristic zero.
\end{remark}

\end{document}